\documentclass{article}
\usepackage[intlimits]{amsmath}
\usepackage{amsfonts}
\usepackage{amssymb}
\usepackage{amsthm}
\usepackage{latexsym}

\newtheorem{Satz}{Theorem}
\newtheorem{Kor}[Satz]{Corollary}
\newtheorem{Lemma}[Satz]{Lemma}
\theoremstyle{definition}
\newtheorem{Def}[Satz]{Definition}

\newcommand{\higgis}{\ensuremath{G_{n,r}}}
\newcommand{\coW}[1]{\ensuremath{\coWord(#1)}}
\newcommand{\QAut}[1]{\ensuremath{\QA(#1)}}

\newcommand{\Hou}[1]{\ensuremath{\Ho(#1)}}
\newcommand{\coCF}{\textit{co\uppercase{CF}}}

\newcommand{\CF}{\textit{\uppercase{CF}}}
\newcommand{\ZZ}{\ensuremath{\mathbb{Z}}}
\DeclareMathOperator{\QA}{QAut}
\DeclareMathOperator{\Ho}{Hou}

\DeclareMathOperator{\coWord}{coW}

\title{Context-Freeness of Higman-Thompson group's co-word problem}
\author{J. Lehnert \and P. Schweitzer}
\begin{document}
\maketitle
\begin{abstract}
The \emph{co-word problem} of a group $G$ generated by a set $X$ is defined as the set of words in $X$ which do not represent $1$ in $G$. We introduce a new method to decide if a permutation group has context-free co-word problem. We use this method to show, that the Higman-Thompson groups, and therefore the Houghton groups, have context-free co-word problem. We also give some examples of groups, that even have an easier co-word problem. We call this property semi-deterministic context-free. The second Houghton group belongs to this class. 
\end{abstract}

\section{Introduction}
Let $G$ be a group with a finite generating set $\mathit{X}$. The \emph{word problem} of $G$ with respect to $\mathit{X}$, denoted by $W(G,X)$, is the set of all words in $(\mathit{X}^\pm)^*$ which represent the identity in $G$.
By \coW{G,X} we denote the complement of $W(G,X)$, i.e. the set of words in $\mathit{X}^{\pm}$ which do not represent $1$ in $G$. \coW{G,X} is called the \emph{co-word problem} of $G$ with respect to $\mathit{X}$. We say the group $G$ has a context-free word problem (resp. co-word problem), or equivalently belongs to the class \CF\ (resp. \coCF), if and only if $W(G,X)$ (resp. \coW{G,X}) is a context-free language. This makes sense, since these questions are independent of the choice of generators, as was shown by Muller and Schupp in \cite{muller} resp. Holt, Rees, R\"over and Thomas in \cite{cocfagain}. As a consequence we can always choose a generating set, which can easily be handled. 

It is a well-known fact, that context-free languages are precisely the languages which can be recognized by a non-deterministic pushdown automaton (NPDA). For a definition and details see e.g. \cite{ginsberg}. The groups of class \CF\ are classified to be the virtual free groups by Muller and Schupp see (\cite{muller}, \cite{dunwoody}). Thereby groups of class \CF\ can already be recognized by deterministic pushdown automata. 

The groups of class \coCF\ were first studied by Holt, R\"over, Rees and Thomas in \cite{cocfagain}. If a deterministic automaton recognizes a language, then there exists an automaton, which recognizes the complement, so the class \CF\ is contained in the class \coCF. Further it is shown, that the class \coCF\ is closed under taking direct products, restricted standard wreath products with \CF-top group, and passing to finitely generated subgroups and finite index overgroups. They also remarked, that there are no further examples of \coCF\ groups known and conjecture that the class \coCF\ is not closed under taking free products. 

The Houghton groups $H_n$, which are defined in Section \ref{basics}, do not arise out of these constructions, so they are indeed new examples. In \cite{roever} R\"over showed, that all Houghton groups are subgroups of the Higman-Thompson groups $G_{n,r}$, which due to Holt and R\"over \cite{coindexed} have indexed co-word problem, a similar but weaker property. There is no group known, which has indexed co-word problem but does not belong to \coCF. The Houghton groups were untruly suspected to be such groups.  In contrary we will show in Section \ref{Higman}:

\begin{Satz}\label{giscocf}
The Higman-Thompson groups $G_{n,r}$ are \coCF.
\end{Satz}
\begin{Kor}\label{Hn}
Let $n\geq 2$. Then $H_n$ is \coCF.
\end{Kor}
Also in Section \ref{basics} we will define the groups \Hou{G}, which are permutation groups similar to the Houghton groups. In fact \Hou{\mathbb{Z}} is isomorphic to $H_2$.

We will show, that \Hou{F_n}\ is \coCF. This result would also follow from the construction in section~\ref{houghton}, but we show a stronger result in section~\ref{leicht}. There is a class of languages in beetween nondeterministic and deterministic context-free languages. This class is recognized by PDA's that may write non-determinis\-ti\-cal\-ly a word onto the stack before the beginning of the computation and then behave totally deterministic. In the case of languages it is easy to see that this class is a true sub-class of the context-free languages, and it seems to be very reasonable that $\coW{H_3}$ is not of that type, but we don't know how to prove it. We will call this class semi-deterministic context-free languages, and obtain: 
\begin{Satz}\label{PFn}
\coW{\Hou{F_n}} is semi-deterministic context-free.
\end{Satz}
This property is also independent of the choice of generators. This implies that finitely generated subgroups of semi-deterministic context-free subgroups inherit this property. A careful analysis of the proofs in \cite{cocfagain} shows, that they easily translate into proofs for semi-deterministic co-word problems. Thus this property is also closed under taking direct products, finite index overgroups and restricted standard wreath products with \CF -top groups.  

In section \ref{sonst} we discuss further groups whose membership to \coCF\ can be established by our method. Throughout this work we will make frequent use of the fact, that group elements are represented by words in the generators. Thereby we will (abusing notation) not distinguish beetween words and group elements.

\section{Definition of Houghton groups and Higman-Thompson groups}\label{basics}
\begin{figure}
\begin{center}
\begin{picture}(300,130)(0,0)
  \put(20,100){\LARGE$*_0^3$}{\scriptsize
  \put(54,46){\line(1,-1){40}}
  \put(46,46){\line(-1,-1){40}}
  \put(50,54){\line(0,1){55}}
  \put(54,46){\circle{3}\hspace*{.3mm}$(0,3)$}
  \put(62,38){\circle{3}\hspace*{.3mm}$(1,3)$}
  \put(70,30){\circle{3}\hspace*{.3mm}$(2,3)$}
  \put(78,22){\circle{3}\hspace*{.3mm}$(3,3)$}
  \put(86,14){\circle{3}\hspace*{.3mm}$(4,3)$}
  \put(46,46){\circle{3}}\put(25,46){$(0,2)$}
  \put(38,38){\circle{3}}\put(17,38){$(1,2)$}
  \put(30,30){\circle{3}}\put(9,30){$(2,2)$}
  \put(22,22){\circle{3}}\put(1,22){$(3,2)$}
  \put(14,14){\circle{3}}\put(-7,14){$(4,2)$}
  \put(50,54){\circle{3}\hspace*{.3mm}$(0,1)$}
  \put(50,65){\circle{3}\hspace*{.3mm}$(1,1)$}
  \put(50,76){\circle{3}\hspace*{.3mm}$(2,1)$}
  \put(50,87){\circle{3}\hspace*{.3mm}$(3,1)$}
  \put(50,98){\circle{3}\hspace*{.3mm}$(4,1)$}}
  \put(200,100){\LARGE$*^3$}{\scriptsize
  \put(250,50){\circle{3}\hspace*{.5mm}$0$}
  \put(250,50){\line(1,-1){44}}
  \put(250,50){\line(-1,-1){44}}
  \put(250,50){\line(0,1){60}}
  \put(250,61){\circle{3}\hspace*{.3mm}$(1,1)$}
  \put(250,72){\circle{3}\hspace*{.3mm}$(2,1)$}
  \put(250,83){\circle{3}\hspace*{.3mm}$(3,1)$}
  \put(250,94){\circle{3}\hspace*{.3mm}$(4,1)$}
  \put(250,105){\circle{3}\hspace*{.3mm}$(5,1)$}
  \put(242,42){\circle{3}}\put(221,42){$(1,2)$}
  \put(234,34){\circle{3}}\put(213,34){$(2,2)$}
  \put(226,26){\circle{3}}\put(205,26){$(3,2)$}
  \put(218,18){\circle{3}}\put(197,18){$(4,2)$}
  \put(210,10){\circle{3}}\put(189,10){$(5,2)$}
  \put(258,42){\circle{3}\hspace*{.3mm}$(1,3)$}
  \put(266,34){\circle{3}\hspace*{.3mm}$(2,3)$}
  \put(274,26){\circle{3}\hspace*{.3mm}$(3,3)$}
  \put(282,18){\circle{3}\hspace*{.3mm}$(4,3)$}
  \put(290,10){\circle{3}\hspace*{.3mm}$(5,3)$}}
  \qbezier(294,40)(280,50)(280,83)\put(280,83){\vector(0,1){2}}
  \qbezier(206,40)(230,50)(230,83)\put(206,40){\vector(-1,-1){2}}
  \qbezier(230,16)(250,40)(270,16)\put(270,16){\vector(1,-1){2}}
  \qbezier(94,40)(80,50)(80,83)\put(80,83){\vector(0,1){2}}
  \qbezier(06,40)(30,50)(30,83)\put(06,40){\vector(-1,-1){2}}
  \qbezier(30,16)(50,40)(70,16)\put(70,16){\vector(1,-1){2}}
  \put(17,57){$x$}\put(47,20){$y$}\put(88,57){$z$}
  \put(217,57){$s_1$}\put(247,20){$s_2$}\put(288,57){$s_3$}

\end{picture}
\end{center}
\caption{Two possible sets of generators of the Houghton group $H_3$. All elements act as shifts, where the arrows define, which rays shifts into each other. While the set $\{x,y,z\}$ is natural for $*^3_0$, the set $\{s_1,s_2,s_3\}$ is natural for $*^3$. They transform into each other as follows: $r=x,\ s=[x,y]y,\ t=z$.  \label{houghtonbild}}
\end{figure}
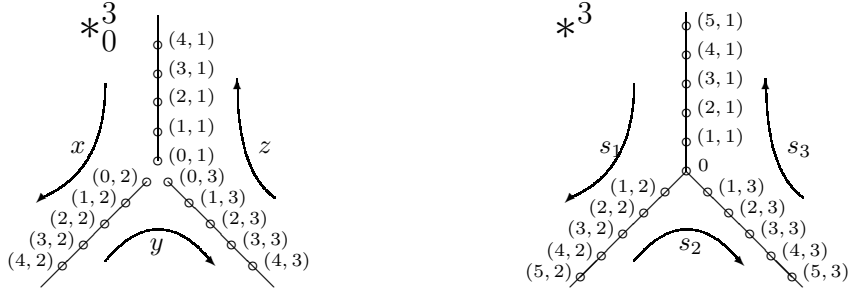

\begin{Def}
Let $\Gamma=(V,E)$ be a locally finite graph. A \textit{quasi-automorphism} of $\Gamma$ is a permutation of $V$, which respects all but finitely many adjacencies. The group of all quasi-automorphisms is denoted by $\QAut{\Gamma}$.
\end{Def}
Note that $\QAut{\Gamma}$ does not change, if we change $\Gamma$ by adding or removing finitely many edges. Thus if $\Gamma$ is a finite graph, \QAut{\Gamma} is just the symmetric group on $V$ which we denote by $S_V$. Let $\mathbb{N}_0$ be the graph with non-negativ integer vertices and $\{i,j\}$ is an edge if and only if $i=j\pm 1$. The disjoint union of $n$ copies of $\mathbb{N}_0$ is called the $n$-star $*^n_0$. We denote the points of $*^n_0$ by a pair $(k,l)$, which is the vertex $k$ of the $l$th copy of $\mathbb{N}_0$ (see figure~\ref{houghtonbild}). For each quasi-automorphism $\Phi$ of $*^n_0$ there exist $r>0,s_1\ldots, s_n\in\mathbb{Z}$ and $\varphi\in S_n$, such that for all $l>r$
\[\Phi\circ(k,l)=(k+s_l,\varphi(l)).\] 
So in the complement of a sufficiently large set the copies of $\mathbb{N}_0$ are permuted and the points are just shifted. Clearly the sum of the $s_i$ has to be $0$.

The $n$th Houghton group $H_n$ is classically defined as the subgroup of all elements $\Phi$ of $\QAut{*^n_0}$ for which this element $\varphi$ is the identity of $S_n$. Since $|S_n|=n!$ we obtain that $H_n$ is a subgroup of finite index in \QAut{*^n_0}. In particular $\QAut{*_0^n}$ is $\coCF$ if and only if $H_n$ is $\coCF$.

As observed obove the group \QAut{*^n_0} does not change if we add edges from the $0$ of the first copy of $\mathbb{N}_0$ to all other $0$'s, and reenumerate $*^n_0$ by adding $1$ to each vertex, which does not belong to the first copy. The result is called $*^n$. For technical reasons, we will use $*^n$ instead of using the usual definition.  Thereby we denote the points of $*^n$ by $0$ and $(k,l),k\in\mathbb{N},l\in\{1,\ldots,n\}$ with the above notation. Be aware, that the shifts in this model, which we define afterwards, are not exactly the same, as in the standard model (see figure~\ref{houghtonbild}). 

Let $n>1$. All $H_n$ are finitely generated (In fact Brown showed in~\cite{brown} that $H_n$ is of type FP$_{n-1}$, but not of type FP$_n$). The shift $s_i\in H_n$ acts on $*_n$ as follows (rays modulo $n$):
\[
s_i(k,l)=\begin{cases}
         (k,l)    & i\neq l\neq i+1\\
         (k-1,l)  & l=i,k\geq 2\\
         0      & l=i,k=1\\ 
         (k+1,l)  & l=i+1
         \end{cases}
\]
and $s_i(0)=(1,i+1)$.  For $n\geq 3$ the commutator $[s_i,s_{i+1}]$ acts as the transposition of $0$ and $(1,i)$. So the group generated by the shifts contains all finite permutations and is thereby equal to $H_n$.

A special case is $H_2$, which needs an additional generator. Because $*^2$ is just the Cayley graph of $\mathbb{Z}$ generated by $X=\{1\}$, we will treat it as a special case of the following definition.
\begin{Def}
Let $G$ be a group finitely generated by $X$. With \Hou{G} we denote the group of those quasi-automorphims of the Cayley graph of $(G,X)$ which act on all but finitely many vertices like a left-multiplication with an element of $G$. In other words $\Hou{G}=S_G\rtimes G$. Here $S_G$ denotes the group of all permutations with finite support.
\end{Def}
Note that this definition is independent of the choice of generators.

If $G$ is generated by $X=\{x_1,\ldots,x_n\}$ then \Hou{G} is generated by $\mathcal{X}=\{x_1,\ldots x_n,\sigma_1,\ldots,\sigma_n\}$. Here $\sigma_i$ is the quasi-automorphism which maps $1$ to $x_i$, $x_i$ to $1$ and stabilizes all other points. Thus $H_2=\Hou{\mathbb{Z}}$ is generated by 2 elements.

From a related concept the \emph{Higman-Thompson groups} arise. Fix two integers $n\geq 2,r\geq 1$ and let $Q=\{q_1\ldots q_r\}$ and $\Sigma=\{\sigma_1\ldots \sigma_n\}$ be two finite sets. Let $\Omega=Q\Sigma^\mathbb{N}$ be the set of infinite sequences starting with some $q_i$ followed by elements of $\Sigma$. A \emph{barrier} is a subset $B\subset Q\Sigma^*$ such that for every element $\omega\in\Omega$ there is exactly one $b\in B$ with $\omega\in b\Sigma^\mathbb{N}$, i.e. $b$ is the unique prefix of $\omega$ in $B$.

Two barriers $B_1,B_2$ of the same cardinality and a bijection $\phi:B_1\rightarrow B_2$ induce a bijection $g_\phi:\Omega\rightarrow\Omega$ by prefix replacement. The set of all such induced bijections forms a group, the Higman-Thompson group $G_{n,r}$.  For a more detailed description of $G_{n,r}$ see \cite{roever}.
 
\section{$\Hou{F_n}$ is semi-deterministic \coCF}\label{leicht}

In order to prove Theorem 2, it will suffice to construct a PDA, which recognizes \coW{\Hou{F_n}}. To demonstrate the way this automaton works, we first solve the co-word problem of $H_2$.
\begin{Satz}
Let $X$ be the generating set of $H_2$ consisting of $X=\{\tau,t\}$, with $\tau=0\leftrightarrow 1$ and $t=(i\mapsto i+1)$.
The Houghton group $H_2$ is semi-deterministic \coCF\ with respect to the generating set $X$.
\end{Satz}
\begin{proof}
We will use the following idea. An element $g\in H_2$ is a pair $(\sigma,s)$ with $\sigma\in S_\infty$ and $s\in\ZZ$ (where we make use of the fact, that $H_2\cong S_\infty\rtimes\ZZ$), hence $g$ is nontrivial if and only if $\sigma$ is a nontrivial permutation or $s\neq 0$.  

The NPDA $P$, which recognizes $CoW(G,X)$, first guesses non-de\-ter\-mi\-nis\-tically whether to check the non-triviality of $s$ or the nontriviality of $\sigma$. In order to check the non-triviality of $s$, all $P$ has to do, is to check that the exponent-sum of $t$ is distinct from $0$. This clearly can be done by a deterministic PDA $P_\ZZ$.

Now we describe, how an automaton $P_S$ decides, if a an element of the form $(\sigma,0)$ is nontrivial. An element $\sigma$ is nontrivial if and only if there is a $k\in\ZZ$ with $\sigma(k)\neq k$. So all $P_S$ has to do, is to guess non-deterministically for which $k$ this condition is tested. 

$P_S$ consists of four states called $I$ (the initial state), $+$, $-$, and $q_A$ (the accept state). Let $\Gamma = \{A,B,C,D,\#\}$ be the stack alphabet, where $\#$ is the end-of-stack symbol. In the beginning, $P_S$ writes a random number of $A$'s (for $k>0$) or $B$'s (for $k<0$) onto the stack. The further behavior of the automaton is completely deterministic. Instead of describing the transition matrix, we explain the main idea. The stack 'remembers' the image of the element $k$, which equals the number of $A$'s (resp. $B$'s) that were written onto the stack in the beginning. So roughly speaking $t$ will add an $A$ and $t^{-1}$ will add a $B$. A $\tau$ only changes the stack, when the stack contains at most one $A$ and no $B$.  

Below these $A$'s and $B$'s the automaton remembers the number of $\tau$ with an effect on the given start-element by $C$'s, when the element is moved to the right and by $D$'s, when its moved to the left. This is possible, because $\tau$ only has an effect on our element, when the stack contains at most one $A$. When the word is read, all the $A$'s and $B$'s are dropped out of the stack and the word is accepted if $C$'s or $D$'s remain. Clearly, during this process $A$ and $B$ cancel each other, so do $C$ and $D$.

The explicit proof can be found in figure~\ref{h2}, which shows the transition matrix of this automaton.
\begin{figure}
\[\begin{array}{c|c|ccc|c}
\textrm{\bf state}&\textrm{\bf read}&\textrm{\bf top of}&\textrm{\bf transforms}&\textrm{\bf state}&\textrm{\bf operation on}\\
&\textrm{\bf letter}&\textrm{\bf stack}&\textrm{\bf to}&&\textrm{\bf top of the stack}\\
\hline
I&t&\Gamma\setminus\{B\}&\leadsto&I&\textrm{write } A\\
I&t&B&\leadsto&I&\textrm{delete }B\\
I&t^{-1}&A&\leadsto&I&\textrm{delete }A\\
I&t&\Sigma\setminus\{A\}&\leadsto&I&\textrm{write } B\\
I&\tau^{\pm 1}&A&\leadsto&+&\textrm{delete } A\\
I&\tau^{\pm 1}&B&\leadsto&I\\
I&\tau^{\pm 1}&C,\Omega&\leadsto&-&\textrm{write } C\\
I&\tau^{\pm 1}&D&\leadsto&-&\textrm{delete } D\\
I&\S&A,B&\leadsto&I&\textrm{delete } A,B\\
I&\S&C,D&\leadsto&q_{A}\\
I&\S&\#& \textrm{reject}&\\
+&\epsilon&A&\leadsto&I&\textrm{write }A\\
+&\epsilon&C&\leadsto&I&\textrm{delete }C\\
+&\epsilon&D,\#&\leadsto&I&\textrm{write }D\\
-&\epsilon&\Gamma&\leadsto&I&\textrm{write }A
\end{array}\] 
\caption{The transition matrix of the automaton, which recognizes \coW{H_2}. $\$$ denotes a symbol, that indicates the end of the input string. Clearly the states $+$ and $-$ are only needed to get account to the second level of the stack. An automaton with a look-ahead of 1 would not need them.\label{h2}} 
\end{figure}
\end{proof}

\begin{proof}[Proof of Theorem~\ref{PFn}]
For $\Hou{F_n}=S_{F_n}\rtimes F_n$ we choose the set of generators $\{x_1, \ldots, x_n, \sigma_1,\ldots,\sigma_{n}\}$. The letters $x_i$ correspond to the generators of the free group and act by (left-) multiplication on the Cayley-graph of $F_n$ and $\sigma_i$  is the switch $1\leftrightarrow x_i$.

The PDA which recognizes \coW{\Hou{F_n}} now does nearly the same as the one that recognizes \coW{H_2}. First it guesses whether it checks the $F_n$ part or the permutation part. The deterministic PDA, which checks the $F_n$ part only has to write the $x_i$-part of the word onto the stack and freely reduce whenever possible.

To check the permutation part the automaton guesses a non-fixpoint and writes it onto the stack. Now the word operates letter by letter. The $\sigma_i$ only have an effect if the stack contains at most one letter. So below the position of the point, the effect of the $\sigma_i$ can be stored on the stack by a reduced word in $\tau_i$. In the end the automaton drops the position and if a $\tau$ remains, the word is non-trivial. Further details are left to the reader.
\end{proof}
\section{Houghton groups are \coCF}\label{houghton}
In order to show that the Houghton groups and later on the Higman-Thompson groups are \coCF\ we need a technical lemma. As C.E. Roever pointed out to us, the result is already proved in \cite{cyclic}, so we suppress the proof 
\begin{Lemma} Let $L$ be a context-free language, then $L^\circ:= \{ yx |
x,y \in \Sigma^\ast, xy \in L\}$, the set of all cyclic
permutations of the words in $L$, is also a context-free language.
\end{Lemma}
Using this we are able to proof that the Houghton groups are \coCF. Despite the fact that this statement follows from Theorem \ref{giscocf}, we carry out the proof explicitly, as it serves as an illustration of the proof of Theorem~\ref{giscocf}.

\begin{proof}[Proof of Corollary~\ref{Hn}]
The Houghton group $H_n$ for $n\geq 3$ is generated by the shifts $s_1,\ldots,s_{n}$ as defined above.
So $s_i$ shifts ray $i$ into ray $(i+1)$. 

Let $\Sigma=\{s_1^{\pm},\ldots,s_{n}^{\pm}\}$,  and $L = \{w\in \Sigma^\ast| w(0) \neq 0\}$ be the set of words for which 0 is not a fixpoint.

 We only need to verify two things now: $L$ is
context free and $L^\circ  = coW(H_n)$, that is every nontrivial
word in the generators is a cyclic permutation
of a word in $L$.

\textbf{Step 1:}  $L^\circ = coW(H_n)$.

\noindent  '$ \subseteq$'  is obvious since conjugation does not
trivialize an element.

\noindent'$\supseteq$'  Let $w$ be a word in $coW(H_n)$. We need to show
there are words $x,y$ with $w=xy$ and a non fixpoint $p$ of $w$, such that
$x(p)= 0$. Let $q=(q_1,q_2)$ be a non-fixpoint. If there is a prefix $x$ of $w$ with $x(q)=0$, we are done. If not, $w$ acts far outside on the ray $(\cdot,q_2)$ like a shift. So there exists also a ray $(\cdot,k)$, on which $w$ acts like an outbound shift (again only far outside). But then there is at least one point $p=(p_1,p_2),p_2\neq k$ whose image is on ray $(\cdot,k)$. This implies the existence of a prefix $x$ of $w$ with $x(p)=0$. 

\textbf{Step 2:} $L$ is context-free. 

\noindent $L$ is even a one counter
language. In one state one can memorize on which ray 0 has
currently been permuted, and on the counter one can memorize how
far.

This completes the proof.
\end{proof}
The group $H_\infty$ is not finitely generated, so it cannot be \coCF. Let $\mathcal{H}_\infty = H_\infty\rtimes\mathbb{Z}$, where $a\in\mathbb{Z}$ acts on $*^\infty$ as the rotation which maps ray $i$ to ray $i+a$. $\mathcal{H}$ is \coCF\ and thereby $H_\infty$ is locally \coCF. The proof is in the same manner and left to the reader. 
\section{The Higman-Thompson groups $\higgis$}\label{Higman}
Using the method of Section \ref{Hn} we are now able to proof Theorem \ref{giscocf}. 
\begin{proof}[Proof of Theorem~\ref{giscocf}]
Higman showed in \cite{higgi} that $G_{n,r}$ is finitely presented. Let $X$ be a finite set of generators for $G_{n,r}$.  For each element $\tau\in \higgis$ there exists an $k_\tau$ such that for all $\omega\in\Omega$ there exist $u,u' \in Q \Sigma ^{\ast}$,  $lg(u)\leq k_\tau$ and $\omega'\in \Sigma^\ast$ with $u\omega'= \omega$ and $\tau(\omega) = u'\omega'$. Here $lg(u)$ denotes the length of $u$. This means, that sequences are only changed on prefixes shorter than $k_\tau$. We observe that $\tau$ induces a map from $Q\Sigma^{k_\tau}\Sigma^*$ (all finite strings with length greater than $k_\tau$) to $Q\Sigma^*$.
Let $k = max_{x\in X} k_x$. 

Let $M = \{qw\sigma_1 \sigma_1 \ldots|q\in Q, w \in \Sigma^{k}\}$ be the set of all infinite sequences for which all entries greater than $k$ are equal to $\sigma_1$. Obviously $|M|=rn^{k}$ and hence  $M$ is finite.

Let $L$ be the set of words in $X^*$, which do not fix all sequences in $M$. 
Again two things remain to be shown.
$L$ is context free and $L^\circ  = coW(\higgis)$: every non-trivial
word in the generators is a cyclic permutation
of a word in $L$.

\textbf{Step 1:} $L^\circ = coW(\higgis)$.

\noindent For every word $z$ in \coW{L} we have to show the existence of a non-fix-sequence which during the action letter-by-letter is mapped onto a sequence of $M$. Obviously $z$ has a non-fix-sequence. Let $r_z=q\sigma_{i_1},\sigma_{i_2}\ldots$ be such a sequence. Let $m\geq k$ be the minimal number such that no prefix of $z$ induces a map which sends the string $q\sigma_{i_1},\ldots\sigma_{i_m}$ to a string of length less than $k$. Because of the minimality of $m$ there exists a prefix $x$ of $z$ which maps $q\sigma_{i_1},\ldots\sigma_{i_m}$ to a string $s$ with length exactly $k$. By the construction of $k$ all sequences with prefix $q\sigma_{i_1},\ldots\sigma_{i_m}$ are mapped onto a sequence with the same prefix as $z(r_z)$ and hence are non-fix-sequences, especially the sequence $\omega=q\sigma_{i_1}\ldots\sigma_{i_m}\sigma_1\sigma_1\ldots$. But $x(\omega)=S\sigma_1\sigma_1\ldots\in M$. This proves step 1.

\textbf{Step 2:}  $L$ is context-free.

\noindent Let $\omega=qw\sigma_1 \sigma_1\ldots $ be an arbitrary element of $M$. Since $M$ is finite, it is sufficient to show, that the language $L'$ of all words in $X^*$, which do not fix $\omega$ is context-free. We describe an automaton with an access depth of $k$ (the automaton is allowed to read and write in the $k$ top levels of the stack), which recognizes $L'$:

The computation starts with the word $qw\sigma_1^l\#$, where $l$ is a random nonnegative integer, on the stack ($q$ on the top). Every letter changes the $k$-prefix in the manner the group-action defines. If the end-of-stack-symbol $\#$ is ever visible to the automaton, it rejects. The word is accepted, if the computation ends with $qw\sigma_1^m\#$ on the stack for some $m$. This can be recognized by the automaton.
\end{proof}
\section{Towards the limits}\label{sonst}
Until now we have not reached the limits of our methods. Because of the statements in the introduction we have also seen, that \QAut{*^n} is \coCF. The arguments of section~\ref{houghton} are clearly strong enough, to show:
\begin{Satz} The following holds:

\begin{enumerate}
\item $\Hou{G}=S_G\rtimes G$ is \coCF\ if and only if G is \coCF. 
\item \QAut{\mathbb{Z}^n} is \coCF.
\end{enumerate}
\end{Satz}
We further conjecture, that \QAut{G} is \coCF\ if and only if $G$ is \coCF\ and \QAut{G} is finitely generated. The motivation for a further study of \QAut{\Gamma} groups comes from the fact, that every finitely generated group is a subgroup of \QAut{\Gamma}-groups. The challenge is just to find a handsome graph $\Gamma$. 

Another reason for a further study of the Higman-Thompson groups is the interesting result of J.C. Birget \cite{birget}, which shows the existence of a subgroup $A\leq G_{3,1}$ and a semidirect product $B=A\rtimes\mathbb{Z}$ such that the co-word problem of $B$ is $\mathop{NP}$-complete and Theorem~\ref{giscocf} implies, that \coW{G_n,r} and the co-word problem of all subgroups can be solved in cubic time.

\section*{Acknowledgements}
We would like to thank S. Rees and C.E. R\"over for posing the questions, which motivated us to work on this subject. Also we would like to thank R. Weidmann for his helpful remarks and discussions.

\end{document}